\documentclass[12pt]{amsart}
\usepackage{amsmath,amssymb,amsthm}
\usepackage[]{latexsym}
\newtheorem{defn0}{Definition}[section]
\newtheorem{prop0}[defn0]{Proposition}
\newtheorem{thm0}[defn0]{Theorem}
\newtheorem{lemma0}[defn0]{Lemma}
\newtheorem{corollary0}[defn0]{Corollary}
\newtheorem{example0}[defn0]{Example}
\newtheorem{remark0}[defn0]{Remark}
\newtheorem{assumption0}[defn0]{Assumption}
\newtheorem{conjecture0}[defn0]{Conjecture}
\newtheorem{notation0}[defn0]{Notation}
\newtheorem{question0}[defn0]{Question}

\newenvironment{definition}{\begin{defn0}\rm}{\end{defn0}}
\newenvironment{proposition}{\begin{prop0}}{\end{prop0}}
\newenvironment{theorem}{\begin{thm0}}{\end{thm0}}

\newenvironment{corollary}{\begin{corollary0}}{\end{corollary0}}

\newenvironment{remark}{\begin{remark0}\rm}{\end{remark0}}

\newenvironment{question}{\begin{question0}\rm}{\end{question0}}

\newcommand{\cA}{{\mathcal A}}
\newcommand{\cQ}{{\mathcal Q}}

\newcommand{\cH}{{\mathcal H}}

\newcommand{\Nm}{\mathrm {Norm}}

\newcommand{\Hom}{{\mathrm {Hom}}}

\newcommand{\PP}{{\mathbb P}}

\newcommand{\cI}{{\mathcal I}}
\newcommand{\cX}{{\mathfrak X}}

\newcommand{\disc}{{\mathrm {disc }}}

\newcommand{\M}{\mathrm{M}}
\newcommand{\Aut}{\mathrm{Aut}}

\newcommand{\End}{{\mathrm{End}}}

\newcommand{\om}{{\omega}}

\newcommand{\Z}{{\mathbb Z}}
\newcommand{\Q}{{\mathbb Q}}
\newcommand{\C}{{\mathbb C}}
\newcommand{\R}{{\mathbb R}}
\newcommand{\F}{{\mathbb F}}

\newcommand{\qbar}{{\bar\Q}}

\newcommand{\osum}{{\oplus}}

\newcommand{\cL}{{\mathcal L}}

\newcommand{\cO}{{\mathcal O}}
\newcommand{\NS}{{\mathrm{NS}}}
\newcommand{\n}{{\mathrm{n}}}
\newcommand{\tr}{{\mathrm{tr}}}

\newcommand{\ra}{{\rightarrow}}

\include{thebibliography}

\begin{document}

\title{Shimura curves embedded in Igusa's threefold}

\author{Victor Rotger}
\footnote{Partially supported by a grant FPI from Ministerio de
Educaci\'{o}n y Ciencia and by Ministerio de Ciencia y Tecnolog\'{\i}a BFM2000-0627}

\address{Universitat Polit\`{e}cnica de Catalunya,
Departament de Matem\`{a}tica Aplicada IV (EUPVG), Av.\
Victor Balaguer s/n, 08800 Vilanova i la Geltr\'{u}, Spain.}

\email{vrotger@mat.upc.es}



\begin{abstract}

Let $\cO $ be a maximal order in a totally indefinite quaternion
algebra over a totally real number field.
In this note we study the locus $\cQ _{\cO }$
of quaternionic multiplication by $\cO $ in the
moduli space $\cA _g$ of principally polarized abelian varieties of
even dimension $g$ with particular emphasis in the two-dimensional case. We
describe $\cQ _{\cO }$ as a
union of Atkin-Lehner quotients of Shimura varieties and we
compute the number of irreducible components of $\cQ _{\cO }$ in terms
of class numbers of CM-fields.

\end{abstract}

\maketitle

\section*{Introduction}

Let $A$ be an abelian variety over a fixed algebraic closure $\qbar $
of the field of rational numbers. By Poincar\'{e}'s Decomposition Theorem, the
algebra of endomorphisms $\End (A)\otimes \Q $ of $A$ decomposes as a direct
sum

$$
\End (A)\otimes \Q \simeq \osum ^s_{i=1} \mathrm{M}_{n_i}(B_i)
$$
of matrix algebras over a division algebra $B_i$.

The ranks $[B_i:\Q ]$ are bounded in terms of
the dimension of $A$ and, by a classical theorem of Albert, the algebras
$B_i$ are isomorphic to either a totally real field, a quaternion algebra over a
totally real field or a division algebra over a CM-field (cf.\,\cite{Mu},
\cite{LaBi}).

We will focus our attention on abelian varieties with totally
indefinite quaternionic multiplication. More precisely, let $F$ be a
totally real number field of degree
$[F:\Q ]=n\geq 1$, let $R_F$ be its ring of integers and
let $\vartheta _F$ denote the different ideal of
$F$ over $\Q $. We also let $F^*_+$ denote the group of totally
positive elements of $F^*$ and $R_{F +}^* = R_F^*\cap F^*_+$.

Let $B$ be a totally indefinite division quaternion algebra over
$F$, i.e. a division algebra of rank $4$ over $F$ such that
$B\otimes _{\Q }\R \simeq \osum _{i=1}^n \M _2(\R )$, and let
$\disc (B)$ denote the reduced discriminant ideal of $B$. We
assume for simplicity that $\vartheta _F$ and $\disc (B)$ are
coprime ideals of $F$. Since $B$ is totally indefinite and
division, it follows from \cite{Vi} that $\disc (B) = \wp _1\cdot
...\cdot \wp _{2 r}$ for pairwise different prime ideals $\wp _i$
of $F$ and $r\geq 1$. We shall denote by $\n = \n _{B/F}$ and $\tr
= \tr _{B/F}$ the reduced norm and reduced trace on $B$, respectively. For
any subset $O$ of $B$, we shall write $O_0 = \{ \beta \in O: \tr
(\beta ) = 0\} $ for the set of pure quaternions of $O$. We shall
also use the notation $O_+ = \{ \beta \in O: \n (\beta )\in
F^*_+\}$ and $O^1 = \{ \beta \in O: \n (\beta )=1\}$.

Finally, let $\cO $ be a maximal order in $B$ and let
$\vartheta (\cO ) = \{ \mu \in \cO : \disc (B)\vert \n (\mu ) \}$ denote
the reduced different of $\cO $ over $R_F$. This is a two-sided ideal
of $\cO $ such that $\n (\vartheta (\cO ))\cdot R_F = \disc (B)$.

Let $\cA _g/\Q $ be the moduli
space of principally polarized abelian varieties of dimension $g = 2 n$. We
let $[(A, \cL )]$ denote the isomorphism class of a principally polarized
abelian variety $(A, \cL )$ regarded as a closed point in $\cA _g$.

It is our aim to investigate the nature of the quaternionic locus

$$
Q_{\cO } \subset \mathcal A _{g}(\C )
$$
of isomorphism classes of complex principally polarized abelian
varieties $[(A, \cL )]$ of dimension $g$ such that $\End (A)\supseteq \cO $.

A reformulation of Proposition 6.2 in \cite{Ro2} yields
that the quaternionic locus $Q_{\cO }$ is {\em nonempty}
if and only if {\em $\disc (B)$ is a totally positive principal ideal of
$F$}. In consequence, for the sake of clarity of
the exposition,
we will assume throughout these notes that the narrow class number of
$F$ is $h_+(F)=1$. This automatically implies the nonemptyness of
$Q_{\cO }$.

{\em Acknowledgements. } I am grateful to A.\ Arenas for some
useful comments and correspondence.

\section{Abelian varieties with quaternionic multiplication}

As a first step in the study of the quaternionic locus
$Q _{\cO }$ in the moduli space $\cA _g$, it is necessary to
understand the geometry of the objects that $Q_{\cO }$ parametrizes.
Let us review some of the results that were accomplished in \cite{Ro2} in
this direction.

\begin{definition}\label{defquat}

An abelian variety with quaternionic multiplication by $\cO $ over $\qbar $ is an
abelian variety $A/\qbar $
such that $\End (A)\simeq \cO $ and $\dim (A) = 2 [F:\Q ] = 2 n$.

\end{definition}

Let $A$ be an abelian variety with quaternionic multiplication by $\cO
$ over $\qbar $ and let $\NS (A)$ denote the free $\Z $-module of
rank $3 n = 3 g/2 $ of line bundles on $A$ up to algebraic equivalence.
We say that two line bundles $\cL $, $\cL'\in \NS (A)$ are isomorphic
if there exists an automorphism $\alpha \in \Aut (A)\simeq \cO ^*$ such
that $\cL ' = \alpha ^*(\cL )$.

\begin{theorem}
\label{NeSe}

Let $A/\qbar $ be an abelian variety with quaternion
multiplication by $\cO $ and let $\iota : \cO \simeq
\End (A)$ be a fixed isomorphism of rings. Then, there is an isomorphism of groups

$$
\begin{array}{clc}
\mathrm{NS}(A)&\stackrel {\sim }{\rightarrow } & \vartheta (\cO )_0 \\
\cL &\mapsto &c_1(\cL ) = \mu
\end{array}
$$
between the N\'{e}ron-Severi group of $A$ and the group of pure
quaternions of the reduced different $\vartheta (\cO )$ of $\cO $.

Moreover, for any two non trivial line bundles $\cL $, $\cL '\in
\NS (A)$, let $\mu =c_1(\cL )$ and $\mu '=c_1(\cL ')$. Then we have that

\begin{enumerate}

\item
$\cL \simeq \cL '$ if and only if there
exists $\alpha \in \cO ^*$ such that $\mu ' = \bar {\alpha }\mu \alpha $.

\item
$\deg (\cL ) = \deg (\varphi _{\cL }:A\ra \hat {A})^{1/2} =
\mathrm{N}_{F/\Q }(\n (\mu )/D )$.

\item
$\cL $ is a polarization if and only if $\n (\mu )\in F^*_+$ and $\mu $
is {\em ample} (cf.\,\cite{Ro2}, \S 5).

\item The Rosati involution on $\cO \stackrel {\iota }{\simeq }\End (A)$ with
respect to $\cL $ is
$$
\begin{matrix}
  \circ : & \cO &\ra & \cO  \\
          & \beta &\mapsto & \mu ^{-1}\bar {\beta }\mu
\end{matrix}
$$

\end{enumerate}
\end{theorem}

Let us write the isomorphism $c_1$ in more explicit terms. Fix an
immersion of $\qbar $ into the field $\C $ of complex numbers and let
$A(\C ) = V/\Lambda $ for some complex vector space $V$ and a lattice $\Lambda
$. Upon the choice of an
isomorphism $\iota : \cO \simeq \End (A)$, the
lattice $\Lambda $ is naturally a left $\cO $-module and, since
$h(F)=h_+(F)=1$, by a theorem of Eichler \cite{Ei55} we know that
$\Lambda \simeq \cO $.

By the Appell-Humbert Theorem, a line bundle $\cL \in \NS (A)$ on $A$
can be regarded as a Riemann form $E: \Lambda \times \Lambda \ra \,\Z $ on
$V$. Let us identify $\Lambda =\cO $ through a fixed isomorphism and
let $t\in F^*$ be any generator of the principal ideal $\vartheta _F$. Then,
the inverse isomorphism $c_1^{-1}: \vartheta (\cO )_{0}
\stackrel {\sim }{\ra }\NS(A)$ maps a pure quaternion $\mu $ to the
Riemann form $E_{\mu }: \cO \times \cO \ra \,\Z $,
$(\beta _1, \beta _2)
\mapsto -\tr _{F/\Q }(\tr (\mu \beta _1 \bar {\beta _2})/\n (\mu ))$.

\begin{remark}

The isomorphism $c_1$ is canonical in the sense that it does not depend
on the choice of any polarization on $A$. However, we warn the reader that it
does depend on the choice of the isomorphism $\iota : \cO \simeq \End
(A)$.

\end{remark}

\begin{theorem}
\label{mainGAV}

Let $A$ be an abelian variety over $\qbar $ with quaternionic multiplication by
$\cO$. Let $D\in F^*_+$ be a totally positive generator of $\disc (B)$.

Then, $A$ is principally polarizable and the
number of isomorphism classes of principal polarizations on $A$ is

$$
\pi _0(A)= \frac {1}{2} \sum _S h(S),
$$
where $S$ runs among the set of orders in the CM-field $F(\sqrt
{-D})$ that contain $R_F[\sqrt {-D}]$ and $h(S)$ denotes its class number.

\end{theorem}

\section{Shimura varieties}\label{Shi}

As in the preceding sections, let $B$ be a
totally indefinite division quaternion algebra over a totally real field $F$
of trivial narrow class number. Let $D\in F^*_+$ be a totally positive
generator of $\disc (B)$.

\begin{definition}

A {\em principally
polarized maximal order} of $B$ is a pair $(\cO , \mu )$ where $\cO \subset B$
is a maximal order and $\mu \in \cO $ is a pure quaternion such that
$\mu ^2 + u D = 0$ for some $u\in R_{F +}^*$.

\end{definition}

Attached to a principally polarized maximal order $(\cO , \mu )$ there
is the following moduli problem over $\Q $: classifying
isomorphism classes of triplets $(A, \iota , \cL )$ given by

\begin{itemize}
\item

An abelian variety $A$ of dimension $g=2n$.

\item

A ring homomorphism $\iota :\cO \hookrightarrow \mbox{ End }(A)$.

\item
A principal pola\-riza\-tion $\cL $ on $A$ such that

$$
\iota (\beta )^{\circ } = \iota (\mu ^{-1}\bar {\beta }\mu )
$$
for any $\beta \in \cO $, where
$\circ :\End (A)\ra \End (A)$ is the Rosati involution with respect to $\cL
$.

\end{itemize}

A triplet $(A, \iota , \cL )$ will be referred to as a polarized abelian variety
with multiplication by $\cO $. Two
triplets $(A_1, \iota _1, \cL
_1)$, $(A_2, \iota _2, \cL _2)$ are isomorphic if there exists an
isomorphism $\alpha \in \Hom (A_1, A_2)$ such that $\alpha \iota
_1 (\beta ) = \iota _2(\beta ) \alpha $ for any $\beta \in \cO $
and $\alpha ^* (\cL _2) = \cL _1\in \NS (A_1)$.
Note also that, since a priori there is no canonical
structure of $R_F$-algebra on $\End (A)$, the immersion $\iota
:\cO \hookrightarrow \End (A)$ is just a homomorphism of rings.

As it was proved by Shimura, the
corresponding moduli functor is coarsely represented by an
irreducible and reduced quasi-projective scheme $\cX _{\mu }/\Q $
over $\Q $ and of dimension $n=[F: \Q ]$.
Moreover, since $B$ is division, the Shimura variety $\cX
_{\mu }$ is complete (cf.\,\cite{Sh63}, \cite{Sh67}).

Let $\mathfrak H = \{ z\in \C : \mathrm{Im }(z)>0 \}$ denote the upper half
plane. Complex analytically, the manifold $\cX _{\mu }(\C )$ can be
described independently of the choice of $\mu $ as the quotient
$$
\cO ^1\backslash \mathfrak H^n\simeq \cX _{\mu }(\C )
$$
of the symmetric space $\mathfrak H^n$ by
the action of the group $\cO ^1$ regarded suitably as a discontinuous
subgroup of $\mbox{SL}_2(\R )^n$. See \cite{Ro3} for details.

In addition, we are also interested in the
Hilbert modular reduced scheme $\cH _F/\Q $ that coarsely represents the functor
attached to the moduli problem of classifying principally
polarized abelian varieties $A$ of dimension $g$ together
with an homomorphism $R_F\hookrightarrow \End (A)$. The Hilbert modular variety
$\cH _F$ has dimension $3 n$ and $\cH _F(\C )$ is the quotient of
$n$ copies of the three-dimensional Siegel half space $\mathfrak H_2$ by a
suitable discontinuous group (cf.\,\cite{Sh63}, \cite{LaBi}).

Notice that, when $F=\Q $, $\cH _F = \cA _2$ is Igusa's
three-fold, the moduli space of principally polarized abelian surfaces.

There are natural morphisms

$$
\begin{matrix}
 \pi : & \cX _{\mu } & \stackrel{\pi _F}{\longrightarrow }&
 \cH _F &\longrightarrow  & \cA _g  \\
       &(A, \iota , \cL )& \mapsto & (A, \iota |_{R_F}, \cL ) &\mapsto & (A, \cL )
\end{matrix}
$$
from the Shimura variety $\cX _{\mu }$ to the Hilbert modular variety $\cH _F$
and the moduli space $\cA _g$ that
consist of gradually {\em forgetting} the quaternionic endomorphism structure. These
morphisms are representable, proper and defined over the field $\Q$ of
rational numbers.

As it will be convenient for our purposes in the rest of this paper, we
introduce the variety $\widetilde {\cX }_{\mu} = \pi (\cX _{\mu })$ to
be the image of the Shimura variety $\cX _{\mu }$ attached to a pure
quaternion $\mu \in \cO $ with $\mu ^2+u D=0$ for some $u\in R_{F +}^*$
in the moduli space $\cA _g$ by the forgetful map $\pi $. It is
important to remark that although
the complex analytical structure of $\cX _{\mu }$ does
not depend on the choice of $\mu $, the construction of the
forgetful map $\pi $ and the subvariety $\widetilde{\cX }_{\mu }$ of $\cA
_g$ do.

The varieties $\widetilde {\cX }_{\mu }$ are
reduced, irreducible, complete and possibly
singular schemes over $\Q $ of dimension $n$. The set of
singularities of $\widetilde {\cX }_{\mu }$ is a finite set
and all the singularities are of {\em quotient type} (cf.\,\cite{vdGe}
for the terminology).

\section{The birational class of the forgetful maps}\label{bir}

It is our aim now to describe the forgetful maps $\pi : \cX _{\mu
} \stackrel{\pi _F}{\longrightarrow }\cH _F \longrightarrow \cA
_g$ as the projection of the Shimura varieties $\cX _{\mu }$ onto
their quotient by suitable Atkin-Lehner groups up to a birational
equivalence. The following groups of Atkin-Lehner involutions were
introduced in \cite{Ro4}. We keep the notations and assumptions of
the Introduction.

\begin{definition}

The {\em Atkin-Lehner group} $W$ of the maximal order $\cO $ is

$$
W = \Nm _{B_+^*}(\cO )/F^*\cdot \cO ^1.
$$

\end{definition}

It was shown in \cite{Ro3} that $W\simeq \Z /2 \Z \times
\stackrel{2 r}{...}\times \Z /2 \Z $, where $2 r$ is the number of
ramified prime ideals of $B$.

\begin{definition}

Let $(\cO , \mu )$ be a principally polarized
maximal order in $B$. A {\em twist} of $(\cO ,\mu )$ is an element
$\chi \in \cO \cap
\Nm _{B^*}(\cO )$ such that $\chi ^2+\n (\chi )=0$ and $\mu \chi =-\chi \mu
$.

\end{definition}

In other words, a twist of $(\cO , \mu )$ is a pure
quaternion $\chi \in \cO \cap \Nm _{B^*}(\cO )$ such that
$$
B = F+F\mu +F\chi +F \mu \chi = (\frac {-u D, -\n (\chi )}{F}).
$$

We say that a principally polarized maximal order $(\cO , \mu )$ in $B$ is
{\em twisting} if it admits some twist $\chi $ in $\cO $. We say
that a maximal order $\cO $ is {\em twisting} if
there exists $\mu \in \cO
$ such that $(\cO , \mu )$ is a twisting principally polarized order.
Finally, we say that $B$ is {\em twisting} if there exists a
twisting maximal order $\cO $ in $B$. Note that $B$ is twisting if and
only if $B\simeq (\frac {-u D, m}{F})$ for some $u\in R^*_{F +}$ and $m\in
F^*$ such that $m\vert D$.

\begin{definition}
\label{twist}

A {\em twisting involution} $\om \in W$ of $(\cO , \mu )$
is an Atkin-Lehner involution such that
$[\om ]=[\chi ]\in
W$ is represented in $B^*$ by a twist $\chi $ of $(\cO , \mu )$.

\end{definition}

We let $V_0=V_0(\cO , \mu )$ denote the subgroup of $W$ generated by the twisting
involutions of $(\cO , \mu )$. For a principally polarized maximal
order $(\cO , \mu )$, let $R_{\mu }=F(\mu )\cap \cO $ and let $\Omega =
\Omega (R_{\mu })=\{ \xi \in R_{\mu }: \xi ^f=1, f\geq 1\} $ denote the
finite group of roots of unity in the CM-quadratic order $R_{\mu }$ over $R_F$.

\begin{definition}

The {\em stable group} of $(\cO , \mu )$ is the subgroup

$$
W_0=U_0\cdot V_0
$$
of $W$ generated by

$$
U_0 = U_0(\cO , \mu )=\Nm _{F(\mu )^*}(\cO )/F^*\cdot \Omega (R_{\mu }),
$$
and the group $V_0$ of twisting involutions of $(\cO , \mu )$.

\end{definition}

As it was also shown in \cite{Ro3}, for any principally polarized maximal order
$(\cO , \mu )$, there are natural monomorphisms of groups
$V_0\subseteq W_0\subseteq W\subseteq \Aut _{\Q }(\cX _{\mu })\subseteq
\Aut _{\qbar }(\cX _{\mu }\otimes \qbar )$.
The question whether the two latter immersions are actually
isomorphisms was studied by the author for the case of Shimura curves in
\cite{Ro1}.

The following was proved in \cite{Ro3}.

\begin{theorem}
\label{mainGSV}

Let $(\cO , \mu )$ be a principally polarized maximal order in $B$ and
let $\cX _{\mu }$ be the Shimura variety attached to it. Then there is
a commutative diagram of finite maps
$$
\begin{matrix}
\cX _{\mu } &  & \stackrel {\pi _F}{\longrightarrow } &  & \cH _F \\
       &  \searrow &                 & \stackrel {b_F}{\nearrow  } &  \\
       &          & \cX _{\mu }/W_0,  &  &
\end{matrix}
$$
where $\cX _{\mu }\rightarrow \cX _{\mu }/W_0$ is the natural projection and
$b_F: \cX _{\mu }/W_0\rightarrow \pi _F(\cX _{\mu })$ is a birational morphism
from $\cX _{\mu }/W_0$ onto the image of $\cX _B$ in $\cH _F$.

The domain of definition of $b_F^{-1}$ is $\pi _F(\cX _{\mu
})\setminus \mathcal T _F$, where $\mathcal T _F$ is a finite set
(of Heegner points).

\end{theorem}

\section{The quaternionic locus}

As above, we let $\cO $ be a maximal order in a
totally indefinite division quaternion algebra $B$ over a totally real field $F$
of trivial narrow class number and we fix a generator $D\in F^*_+$ of $\disc (B)$.

It is the aim of this section to use the preceding results to study
the quaternionic locus $Q_{\cO
}$ in $\cA _g(\C )$. As we mentioned in the Introduction, since
$h_+(F)=1$, the set $Q_{\cO }$ is not empty.

\begin{definition}\label{Hee}

A Heegner point in $Q_{\cO }$ is an isomorphism class $[(A, \cL )]$ of
a principally polarized abelian variety such that $\End (A)\varsupsetneq \cO
$.

\end{definition}

According to the Definitions \ref{defquat} and \ref{Hee}, we note
that $Q_{\cO }$ is the disjoint union of the set of principally
polarized abelian varieties $[(A, \cL )]$ with quaternionic
multiplication by $\cO $ and the set of Heegner points. The former
is a subset of $Q_{\cO }$ whose closure with respect to the
analytical topology is $Q_{\cO }$ itself. The latter is also a dense but
discrete subset of $Q_{\cO }$ (cf.\,\cite{Sh67}).

In order to understand the nature of the locus $Q_{\cO }$, we observe that
for any principally polarized pair $(\cO , \mu )$, the
set $\widetilde {\cX }_{\mu }(\C )$ of complex points of
the Shimura variety $\widetilde {\cX }_{\mu }/\Q $ attached to $(\cO , \mu )$
sits inside $Q_{\cO }$.

\begin{proposition}\label{inter}

Let $\mu $, $\mu '\in \cO _0$ be two pure quaternions such that
$\n (\mu ) = u D$ and $\n (\mu ') = u' D$ for some units $u$, $u'\in
R_{F +}^*$. If $\widetilde {\cX }_{\mu }(\C )$ and
$\widetilde {\cX }_{\mu '}(\C )$ are different subvarieties of $\cA _{g}(\C
)$, then $\widetilde {\cX }_{\mu }(\C )\cap \widetilde {\cX }_{\mu
'}(\C )$ is a finite set of Heegner points.

\end{proposition}

{\em Proof. } Assume that the isomorphism class $[A, \cL ]$
of a principally polarized abelian variety falls at the intersection of
$\widetilde {\cX }_{\mu }$ and $\widetilde {\cX }_{\mu '}$ in $\cA _g$. Write
$[A, \cL ]=\pi ([A, \iota ,\cL ]) = \pi ([A', \iota ',\cL '])$ as the
image by $\pi $ of points in $\cX _{\mu }$ and $\cX _{\mu '}$,
respectively. Since $[(A, \cL )]=[(A', \cL ')]\in
\cQ _{\cO }$, we can identify the pair $(A, \cL ) = (A', \cL ')$
through a fixed isomorphism of polarized abelian varieties.

Let us assume that
$[A, \cL ]=[A', \cL ']$ was not a Heegner
point. Then $\iota :\cO \simeq \End (A)$ would be an isomorphism of rings
such that $c_1(\cL ) = \mu $. We then would have
by Theorem \ref{NeSe}, \S 4, that $c_1(\cL ) = c_1(\cL ') = \mu = \mu' $
up to multiplication by elements in $F^*$. Since
$\widetilde {\cX }_{\mu } = \widetilde {\cX }_{u \mu }$ for all units $u\in R_F^*$,
this would contradict the statement. Since
the set of Heegner points in $\cA _g(\C )$ is discrete, we conclude that
any two irreducible components of $\cQ _{\cO }$ meet at a finite set of
Heegner points. $\Box $

\begin{proposition}
\label{locus}

\begin{enumerate}

\item The locus $Q_{\cO }$ is the set of complex points $\mathcal Q
_{\cO }(\C )$ of a reduced complete subscheme $\mathcal Q _{\cO }$
of $\mathcal A_g$ defined over $\Q $.

\item Let $\rho (\cO )$ be the number of
absolutely irreducible components of $\mathcal Q _{\cO }$. Then
there exist quaternions $\mu _k\in \cO _0$ with
$\mu _k^2 + u_k D=0$ for $u_k\in R_{F +}^*$, $1\leq k\leq \rho (\cO )$,
such that

$$
\mathcal Q _{\cO } = \bigcup
\widetilde {\cX }_{\mu _k}.
$$
is the decomposition of $\cQ _{\cO }$ into irreducible components.

\end{enumerate}

\end{proposition}

{\em Proof. } Let $[(A, \cL )]\in Q_{\cO }$ be the isomorphism class of
a complex principally polarized abelian variety such that $\End (A)\simeq \cO
$. Fix an isomorphism $\iota :\cO \simeq \End (A)$.
By Theorem \ref{NeSe}, \S 4, the
Rosati involution with respect to $\cL $ on $\cO $ must be of
the form $\varrho _{\mu }:\cO \ra \cO $, $\beta \mapsto \mu ^{-1}\bar {\beta }\mu $
for some $\mu \in \cO $ with $\mu ^2 + u D
=0$, $u\in R_{F +}^*$. Thus $[(A, \cL )] = \pi ([A, \iota , \cL ])\in
\widetilde {\cX }_{\mu }(\C )$, namely the set
of complex points on a reduced, irreducible, complete and possibly singular scheme
over $\Q $ (cf.\,\cite{Sh63}, \cite{Sh67}). Since the set of Heegner
points $[(A, \cL )] \in \widetilde {\cX }_{\mu }(\C )$
is a discrete set which lies on the Zariski
closure of its complement, we conclude that $Q_{\cO }$ is the union
of the Shimura varieties $\widetilde {\cX }_{\mu }(\C )$ as $\mu $
varies among pure quaternions satisfying the above properties.

Let us now show that $Q_{\cO }$ is actually covered by
finitely many pairwise different Shimura varieties. Let
$A/\C $ be an arbitrary abelian variety with quaternionic
multiplication by $\cO $ and fix an isomorphism $\iota : \cO \simeq \End (A)$.

Let $(\cO , \mu )$ be any principally polarized pair. Since $h_+(F)=1$, there
exists a unit $u\in R_F^*$ such that $u \mu $ is an ample quaternion in the
sense of \cite{Ro2}, \S 5. Let $\cL \in \NS (A)$ be
the line bundle on $A$ such that $c_1(\cL )^{-1}= u \mu $. From
Theorem \ref{NeSe}, it
follows that $\cL $ is a principal
polarization on $A$ such that the isomorphism class of the triplet
$(A, \iota , \cL )$ corresponds to a closed point in $\cX _{\mu }(\C
)$ and hence $[A, \cL ]\in \widetilde {\cX }_{\mu }$.

Since, by Proposition \ref{inter},
the intersection points of two different Shimura varieties
$\widetilde {\cX }_{\mu }(\C )$ and $\widetilde {\cX }_{\mu '}(\C )$ in
$\cA _{g }(\C )$ are Heegner points, this
shows that for every irreducible component of $Q _{\cO }$ there exists at least
one principal polarization $\cL $ on $A$ such that $[A, \cL ]$ lies on
it. Consequently, the number $\pi _0(A)$ of isomorphism classes of
principal polarizations on $A$ is an upper bound for the number $\rho (\cO
)$ of irreducible components of $Q_{\cO }$. Since, by Theorem \ref{mainGAV}, the
number $\pi _0 (A)$ is a finite number, this yields the
proof of the proposition. $\Box $

$\\ $
In view of Proposition \ref{locus}, it is natural to pose the following

\begin{question}\label{Qirr}

What is the number $\rho (\cO )$ of irreducible components of $\cQ _{\cO
}$? When is $\cQ _{\cO }$ irreducible?

\end{question}

\section{The distribution of principal polarizations on an abelian variety
in $\cQ _{\cO }$}

Let us relate Question \ref{Qirr} to the following problem.
In Theorem \ref{mainGAV}, we computed the number $\pi _0(A)$ of
principal polarizations on an abelian variety $A$ with quaternion
multiplication by $\cO $ as the finite sum of relative class numbers of
suitable orders in the CM-fields $F(\sqrt {- u D})$ for $u\in
R_{F +}^*/R_F^{*2}$. This has the following modular interpretation:

Let $\cL _1$, ..., $\cL _{\pi _0(A)}$ be representatives of the $\pi _0(A)$
distinct principal polarizations on $A$. Then
the pairwise nonisomorphic principally polarized abelian varieties
$[(A, \cL _1)]$, ..., $[(A, \cL _{\pi _0(A)})]$ correspond to all closed points
in $\cQ _{\cO }$ whose underlying abelian variety is isomorphic to
$A$. We then naturally ask the following

\begin{question}\label{Qdistri}

Let $A$ be an abelian variety with quaternionic multiplication by
$\cO $. How are the distinct principal polarizations $[(A, \cL _j)]$ distributed
among the irreducible components $\widetilde {\cX }_{\mu _k}$ of $\cQ _{\cO }$?

\end{question}

It turns out that the two questions above are related. The
linking ingredient is provided by the definition below, which
establishes a slightly coarser equivalence relationship on polarizations
than the one considered in Theorem \ref{NeSe}, $\S 1$.

\begin{definition}

Let $A$ be an abelian variety with quaternionic multiplication by $\cO
$ over $\qbar $.

\begin{enumerate}
\item
Two polarizations $\cL $ and $\cL '$ on $A$ are weakly isomorphic if
$c_1(\cL ) \simeq m c_1(\cL ')\in \NS (A)$ for some $m\in F^*_+$. We shall denote it
$\cL \simeq _{w}\cL '$.

\item
Two principal polarizations $\cL $ and $\cL '$ on $A$ are
{\em Atkin-Lehner isogenous}, denoted $\cL \sim \cL '$, if there
is an isogeny $\om \in \cO \cap \mathrm{Norm}_{B_+^*}(\cO )$ of $A$ such that

$$
\om ^*(\cL )\simeq _{w}\cL '.
$$

\end{enumerate}
\end{definition}

We note that there is a closed relationship between the
above definition and the modular interpretation of the Atkin-Lehner group $W$
given in \cite{Ro3}.

\begin{definition}

Let $A$ be an abelian variety with quaternionic multiplication by $\cO
$ over $\qbar $. We let $\hat {\Pi }_0(A)$ be the
set of principal polarizations on $A$ up to Atkin-Lehner isogeny and we
let $\hat {\pi }_0(A) = \sharp \hat {\Pi }_0(A)$ denote its
cardinality.

\end{definition}

\begin{theorem} [Distribution of principal polarizations]\label{distri}

Let $A$ be an abelian variety with quaternionic multiplication by
$\cO $ over $\qbar $ and let $\cL _1$, ..., $\cL _{\pi _0(A)}$ be representatives of
the $\pi _0(A)$
distinct principal polarizations on $A$.

Then, two closed points $[A, \cL _i]$ and
$[A, \cL _j]$ lie on the same irreducible component of $\cQ _{\cO }$
if and only if the
polarizations $\cL _i$ and $\cL _j$ are Atkin-Lehner isogenous.

\end{theorem}

{\em Proof. } We know from Proposition \ref{locus} that any irreducible
component of $\cQ _{\cO }$ is $\widetilde {\cX }_{\mu }$ for some
principally polarized pair $(\cO , \mu )$. We single out and
fix one of these.

Let $\cL $ be a principal polarization on $A$ such that $[(A, \cL
)]$ lies on $\widetilde {\cX }_{\mu }$ and let
$\cL '$ be a second principal polarization on $A$. We claim that
$[A, \cL ']\in \widetilde {\cX }_{\mu }$ if and only if there exists
$\om \in \cO \cap \mathrm{Norm}_{B^*_+}(\cO )$ such that $\cL '$ and
$\om ^*(\cL )$ are weakly isomorphic.

Assume first that $\cL '\simeq _{w}\om ^*(\cL )$ for some
$\om \in \cO \cap \mathrm{Norm} _{B^*_+}(\cO )$. This amounts to saying
that $\bar {\om } c_1 (\cL ) \om = m c_1(\cL ')$ for some $m\in F^*$.
Since both $\om ^*(\cL )$ and $\cL '$ are polarizations, we deduce from
Theorem \ref{NeSe}, \S 3, that $m\in F^*_+$. Moreover, since
$\cL $ and $\cL '$ are principal, we obtain from
Theorem \ref{NeSe}, $\S 2$, that $m = u \n (\om )$ for some $u\in
R_{F +}^*$.

Note that $(A, \iota _{\om }, \cL ')$ is a principally polarized
abelian variety with
quaternionic multiplication such that the
Rosati involution that $\cL '$ induces on $\iota _{\om }(\cO )$
is $\varrho _{\mu }$. Indeed, this follows
because $\iota _{\om }(\beta )^{\circ _{\cL '}}=
\iota ((\om ^{-1}\beta \om ))^{\circ _{\cL '}} =
\iota ((\om ^{-1}\mu \om )^{-1}\overline {\om ^{-1} \-\beta \om }\-(\om ^{-1}\mu \om
)) \-= \iota _{\om }(\mu ^{-1}\bar {\beta } \mu )$. This shows that,
if $\cL '$ and $\om ^*(\cL )$ are weakly
isomorphic for some $\om \in \cO \cap \Nm _{B^*_+}(\cO )$, then $[A, \cL ']\in
\widetilde {\cX }_{\mu }$.

Conversely, let us assume that $[A, \cL ']\in \widetilde {\cX }_{\mu }$. Let
$\iota ':\cO \hookrightarrow \End (A)$ be
such that $[A, \iota ', \cL ']\in \cX _{(\cO , \cI _{\vartheta
}, \varrho _{\mu })}$. By the Skolem-Noether
Theorem, it holds that $\iota ' = \om ^{-1}\iota \om $ for some
$\om \in \mathrm{Norm}_{B^*}(\cO )$;
we can assume that $\om \in \cO
$ by suitably scaling it. Since it holds that
$\iota _{\om }(\beta )^{\circ _{\cL '}}= \iota _{\om }(\mu ^{-1}\bar \beta
\mu)$ for any $\beta \in \cO $, we have that
$c_1(\cL ') = u \om ^{-1}c_1(\cL )\om
$ for some $u\in R_F^*$ such that $u \n (\om )\in F^*_+$. Since
$\n (\cO ^*)=R_F^*$, we
can find $\alpha \in \cO ^*$ with reduced norm $\n (\alpha )=u^{-1}$ and thus
$\om \alpha \in B^*_+$.

Let $\cL _{\om \alpha }$ be the polarization on $A$ such that
$c_1(\cL _{\om \alpha }) = \frac {u}{\n (\om )} c_1( (\om \alpha )^*(\cL
))$.  The automorphism $\alpha \in \cO ^*=\Aut (A)$
induces an isomorphism between the polarizations $\cL _{\om \alpha }$
and $\cL '$, since $c_1(\alpha ^*(\cL ')) = \bar {\alpha }(u \om ^{-1}c_1(\cL )\om )
\alpha =c_1(\cL _{\om \alpha })$. Hence $\cL '$ is weakly isomorphic
to $(\om \alpha )^*(\cL )$. This concludes our claim above and
also proves the theorem. $\Box $

\begin{corollary}
\label{genus}

The number of irreducible components of $\cQ _{\cO }$ is

$$
\rho (\cO ) = \hat \pi_0(A),
$$
independently of the choice of $A$.

\end{corollary}

For any irreducible
component $\widetilde {\cX }_{\mu _k}$ of $\mathcal Q _{\cO
}$, let $\Pi _0^{(k)}(A)\subset \Pi _0(A)$ denote the set of
isomorphism classes of the isogeny class of principal polarizations
lying on $\widetilde {\cX }_{\mu _k}$.

As another immediate consequence of Theorem \ref{distri}, the
following corollary establishes an internal structure
on the set $\Pi _0(A)$. Roughly, it asserts that
$\Pi _0(A) = \bigcup _{k=1}^{\rho (\cO )} \Pi _0^{(k)}(A)$ is
the disjoint union of the sets $\Pi _0^{(k)}(A)$, which are equipped with a free
and transitive action of a 2-torsion finite abelian group.

\begin{corollary}

Let $A$ be an abelian variety with quaternionic multiplication by
$\cO $. Let $\widetilde {\cX }_{\mu _k}$ be an irreducible component of
$\cQ _{\cO }$ and let $W_0^{(k)}\subseteq W$ be the stable subgroup
attached to the polarized order $(\cO , \mu _k)$.

Then there is a free and
transitive action of $W/W_0^{(k)}$ on the Atkin-Lehner isogeny class
$\Pi _0^{(k)}(A)$ of principal polarizations lying on $\widetilde {\cX }_{\mu _k}$.

\end{corollary}

In the case of a non twisting maximal order $\cO $, we have that the
stable group $W_0 (\cO , \mu )$ attached to a principally polarized
pair $(\cO , \mu )$ is $U_0(\cO , \mu )$. The following corollary
follows from the proof of Lemma 2.10 in \cite{Ro4}.

\begin{corollary}\label{equi}

Let $\cO $ be a non twisting maximal order in $B$ and assume that, for any
$u\in R_{F +}^*$, any
primitive root of unity of odd order in the CM-field $F(\sqrt {-u D})$
is contained in the order $R_F[\sqrt {-u D}]$.

Let $A$ be an abelian variety with quaternionic multiplication by $\cO $.
Then the distinct isomorphism
classes of principally
polarized abelian varieties $[(A, \cL _1)]$, ..., $[(A, \cL _{\pi _0(A)})]$ are
{\em equidistributed} among the $\rho (\cO )$ irreducible components
of $\mathcal Q _{\cO }$.

In particular, it then holds that

$$
\pi _0 (A) = \frac {|W|}{|W_0|}\cdot \rho (\cO ).
$$

\end{corollary}

\section{Shimura curves embedded in Igusa's threefold}
\label{curves}

The whole picture becomes particularly neat when we consider the simplest
case of quaternion algebras over $\Q $. Let then $B$ be an
indefinite quaternion algebra over $\Q $ of discriminant
$D = p_1\cdot ...\cdot p_{2 r}$ and let $\cO $ be a maximal
order in $B$. Since
$h(\Q )=1$, there is a single choice of $\cO $ up to conjugation by
$B^*$. Moreover, all left ideals of $\cO $ are principal and hence isomorphic
to $\cO $ as left $\cO $-modules.

Let $A$ be a complex abelian surface with quaternionic multiplication
by $\cO $. By Theorem \ref{mainGAV}, $A$ is principally polarizable and
the number of isomorphism classes of principal polarizations on $A$ is

$$
\pi _0(A)= \dfrac {h^{\widetilde {}}(-D)}{2},
$$
where, for any nonzero squarefree integer $d$, we write
$$
h^{\tilde {}}(d) = \begin{cases}
    h(4 d)+h(d)  & \text{ if }d\equiv 1\quad \mbox{ mod }4, \\
    h(4 d) & \text{otherwise.}
  \end{cases}
$$

For any integral element $\mu \in \cO $ such that $\mu ^2 + D = 0$, let
now $\cX _{\mu }$ be the Shimura curve
that coarsely represents the functor which classifies
principally polarized abelian surfaces $(A, \iota , \cL )$ with
quaternionic multiplication by $\cO $
such that the Rosati involution with respect to $\cL $ on
$\cO $ is $\varrho _{\mu }$. This is an algebraic curve over $\Q
$ whose isomorphism class does not actually depend on the choice of the
quaternion
$\mu $, but only on the discriminant $D$ (cf.\,\cite{Sh67}). Hence, it is usual to
simply denote this isomorphism class as $\cX _D$.

Let $W = \{ \omega _m: m\vert D\} \simeq (\Z /2\Z )^{2 r}$
be the Atkin-Lehner group attached to $\cO $ in Section \ref{bir}.
We know that $W\subseteq \Aut _{\Q }(\cX
_D)$ is a subgroup of the group of automorphisms of the Shimura curve $\cX
_D$.

Let now $\cA _2$ be Igusa's three-fold, the moduli space of principally polarized
abelian surfaces. By the work of Igusa (cf.\,\cite{Ig}), it is
an affine scheme over $\Q $ that
contains, as a Zariski open and dense subset, the moduli space $\mathcal M
_2$ of curves of genus $2$, immersed in $\cA _2$ via the Torelli
embedding.

Sitting in $\cA _2$ there is the quaternionic locus $\cQ _{\cO }$ of
isomorphism classes of principally polarized abelian surfaces $[(A, \cL
)]$ such that $\End (A)\supseteq \cO $. Since all maximal orders $\cO $
in $B$ are pairwise conjugate, the quaternionic locus $\cQ _{\cO }$
does not actually depend on the choice of $\cO $ and we may simply
denote it by $\cQ _{\cO } = \cQ _D$.

As explained in Section \ref{Shi} and \ref{bir}, there are
forgetful finite morphisms $\pi :\cX _{\mu }\rightarrow \cQ _D \subset
\cA _2$ which map the Shimura curve $\cX _{\mu }$ onto an irreducible
component $\widetilde {\cX }_{\mu }$ of $\cQ _D$. We insist on the fact
that the image $\widetilde {\cX }_{\mu }\subset \cQ _D$ {\em does depend} on
the choice of the quaternion $\mu $.

Let us now compare the non twisting and twisting case,
respectively. We first assume that

$$
B\not \simeq (\frac {-D, m}{\Q })
$$
for all positive divisors $m\vert D$ of $D$. Then all principally polarized
pairs $(\cO , \mu )$ in $B$ are {\em non twisting} and the stable
subgroup attached to $(\cO , \mu )$ is
$$
W_0 = U_0 = \langle \om _D\rangle \subset W,
$$
independently
of the choice of $\mu $. By Theorem \ref{mainGSV}, we
deduce that any irreducible component $\widetilde {\cX }_{\mu }$ of
$\cQ _D$ is birationally equivalent to the Atkin-Lehner quotient
$\cX _D/\langle \om _D \rangle $ and thus the
quaternionic locus $\cQ _D$ in $\cA _2$
is the union of pairwise birationally equivalent
Shimura curves $\widetilde {\cX }_{\mu _1}$,
...,$\widetilde {\cX }_{\mu _{\rho (\cO )}}$, meeting at a finite set of
Heegner points.

Moreover, for any abelian surface $A$ with quaternionic multiplication
by $\cO $, it follows from
Theorem \ref{distri} that the closed points $\{ [(A, \cL _j)]\} _{j=1}^{\pi
_0(A)}$ are equidistributed among the $\rho (\cO )$ irreducible components
of $\cQ _D$. In addition, Corollary \ref{equi} ensures that

$$
|W/W_0| = 2^{2 r-1} \vert \pi _0(A),
$$
as genus theory for binary quadratic forms already
predicts.

Finally, we obtain that the number of irreducible components of $\cQ _D$
in the non twisting case is

$$
\rho (\cO ) = \dfrac {h^{\widetilde {}}(-D)}{2^{2 r}}.
$$

On the other hand, let us assume that

$$
B\simeq (\frac {-D, m}{\Q })
$$
for some $m\vert D$. In this case, there can be two different
birational classes of irreducible components on $\cQ _D$. Indeed, the assumption
means that there
exist pure quaternions $\mu \in \cO $, $\mu ^2+D=0$, such that $(\cO , \mu
)$ is a {\em twisting} principally polarized maximal order. Then
$$
W_0 (\cO , \mu ) = \langle \om _m, \om _D\rangle
$$
and $\widetilde {\cX }_{\mu }$
is birationally equivalent to $\cX _D/\langle \om _m, \om _D\rangle $.
We may refer to $\widetilde {\cX }_{\mu }$ as a
twisting irreducible component of $\cQ _D$.

In addition to these, there may exist non twisting
polarized orders $(\cO , \mu )$ such that the
corresponding irreducible components $\widetilde {\cX }_{\mu }$ of $\cQ _D$
are birationally equivalent to $\cX _D/\langle \om _D\rangle $. We
may refer to these as the non twisting irreducible components of $\cQ
_D$.

We then have the following lower and upper bounds for the number of
irreducible components of $\cQ _D$:

$$
\dfrac {h^{\widetilde {}}(-D)}{2^{2 r}} <\quad \rho (\cO ) \quad \leq
\dfrac {h^{\widetilde {}}(-D)}{2^{2 r -1}}.
$$

Summing up, we obtain the following

\begin{theorem}\label{pq}

Let $B$ be an indefinite division quaternion algebra over $\Q $
of discriminant $D = p_1\cdot...\cdot p_{2 r}$. Then, the
quaternionic locus $\cQ _D$ in $\cA _2$ is irreducible if and only if

$$
h^{\widetilde {}}(-D) =
  \begin{cases}
    2^{2 r -1} & \mbox{ if }B\simeq (\frac {-D, m}{\Q })\mbox{
    for some }m\vert D, \\
    2^{2 r} & \text{otherwise}.
  \end{cases}
$$

\end{theorem}

{\em Proof. } If $B$ is not a twisting quaternion algebra, we already
know from the above that the number of irreducible components of $\cQ
_D$ is $\dfrac {h^{\widetilde {}}(-D)}{2^{2 r}}$. Hence, in this case, the
quaternionic locus of discriminant $D$ in $\cA _2$ is irreducible if
and only if $h^{\widetilde {}}(-D) = 2^{2 r}$. If on the other hand $B$ is
twisting, it follows from the above inequalities that $\cQ _D$ is
irreducible if and only if $h^{\widetilde {}}(-D) = 2^{2 r - 1}$. $\Box $

In view of Theorem \ref{pq}, there arises a closed relationship
between the irreductibility of the quaternionic locus in Igusa's
threefold and the genus theory of integral binary quadratic forms
and the classical {\em numeri idonei} studied by Euler, Schinzel
and others. We refer the reader to \cite{Ar95} and \cite{Sch59}
for the latter.

\section{Hashimoto-Murabayashi's families}

As the simplest examples to be considered, let $B_6$ and $B_{10}$ be the
rational quaternion algebras of discriminant $D = 2\cdot 3 = 6$ and $2\cdot 5 =
10$, respectively. Hashimoto and Murabayashi \cite{HaMu} exhibited
two families of principally polarized abelian surfaces with
quaternionic multiplication by a maximal order in these quaternion
algebras. Namely, let

$$
C_{(s, t)}^{(6)}: Y^2 = X (X^4 + P X^3 + Q X^2 + R X + 1)
$$
be the family of curves with

$$
P = 2 s + 2 t, \quad Q = \frac {(1+2 t^2) (11-28 t^2 + 8 t^4)}{3 (1-t^2)(1-4
t^2)}, \quad R = 2s -2t
$$
over the base curve

$$
g^{(6)}(t, s) = 4 s^2 t^2 - s^2 + t^2 + 2 = 0.
$$

And let

$$
C_{(s, t)}^{(10)} : Y^2 = X (P^2 X^4 + P^2 (1+R) X^3 + P Q X^2 + P (1-R) X + 1)
$$
be the family of curves with

$$
P = \frac {4 (2 t + 1)(t^2-t-1)}{(t-1)^2},
\quad Q = \frac {(t^2+1)(t^4+8t^3-10t^2-8t+1)}{t (t-1)^2 (t+1)^2}
$$
and
$$
R = \frac {(t-1) s}{t (t+1) (2 t + 1)}
$$
over the base curve

$$
g^{(10)}(t, s) = s^2 - t (t-2) (2 t +1) = 0.
$$

Let $J_{(s, t)}^{(6)} =
\mbox{Jac} ( C_{(s, t)}^{(6)})$ and $J_{(s, t)}^{(10)} =
\mbox{Jac} ( C_{(s, t)}^{(10)})$ be the Jacobian surfaces of the
fibres of the families of curves above respectively. It was proved in \cite{HaMu}
that their ring of endomorphisms contain a maximal order in $B_6$ and
$B_{10}$, respectively.

Both $B_6 = (\frac {-6, 2}{\Q })$ and
$B_{10} = (\frac {-10, 2}{\Q })$ are twisting quaternion algebras. Moreover,
it turns out from our formula for $\pi _0(A)$ in the above section
that any abelian surface $A$ with
quaternionic multiplication by a maximal order in either $B_6$ or
$B_{10}$ admits a {\em single} isomorphism class of principal
polarizations. This implies that
$\rho (B_6)=\hat \pi _0(A) = \pi _0(A) = 1$ and $\rho (B_{10})=\hat \pi _0(A) =
\pi _0(A) = 1$, respectively.

Moreover, the Shimura
curves $\cX _6/\Q $ and $\cX _{10}/\Q $ have genus $0$, although they are not
isomorphic to $\PP ^1_{\Q }$ because there are no rational points on
them. However, it is easily seen that $\cX _6/W_0 = \cX _6/W \simeq \PP ^1_{\Q
}$ and $\cX _{10}/W_0 = \cX _{10}/W \simeq \PP ^1_{\Q }$, respectively.

As it is observed in \cite{HaMu}, the base curves $g^{(6)}$
and $g^{(10)}$ are curves of genus $1$ and not of genus $0$ as it
should be expected. This is explained by the fact that
there are obvious isomorphisms between
the fibres of the families $C^{(6)}$ and $C^{(10)}$, respectively.

Ibukiyama, Katsura and Oort \cite{IbKaOo} proved that the
supersingular locus in $\cA _2/\bar {\F _p}$ is irreducible if and
only if $p\leq 11$. As a corollary to their work, Hashimoto and
Murabayashi obtained that the reduction mod $3$ and $5$ of the
family of Jacobian surfaces with quaternionic multiplication by
$B_6$ and $B_{10}$ respectively yield the single irreducible
component of the supersingular locus in $\cA _2/\bar {\F _3}$ and
$\cA _2/\bar {\F _5}$ respectively. The following statement may be
considered as a lift to characteristic $0$ of these results.

\begin{theorem}\label{uni}

\begin{enumerate}
\item
The quaternionic locus $\cQ _6$ in $\cA _2/\Q $ is absolutely
irreducible and birationally equivalent to $\PP ^1_{\Q }$ over $\Q $.
A universal family over $\qbar $ is given by
Hashimoto-Murabayashi's family $C^{(6)}$.

\item
The quaternionic locus $\cQ _{10}$ in $\cA _2/\Q $ is absolutely
irreducible and birationally equivalent to $\PP ^1_{\Q }$ over $\Q $.
A universal family over $\qbar $ is given by
Hashimoto-Murabayashi's family $C^{(10)}$.

\end{enumerate}
\end{theorem}

{\em Proof. } This follows from Theorem \ref{pq} and the
discussion above. $\Box $

In particular, we obtain from Theorem \ref{uni} that every principally
polarized abelian surface $(A, \cL)$ over $\qbar $ with quaternionic multiplication
by a maximal order of discriminant $6$ or $10$ is isomorphic over $\qbar $ to
the Jacobian variety of one of the curves $C_{(s, t)}^{(6)}$ or $C_{(s,
t)}^{(10)}$, except for finitely many degenerate cases.

\end{document}